\documentclass[a4paper,10pt,reqno]{amsart}
\usepackage{verbatim}       
\usepackage{amssymb}       
\usepackage{mathrsfs}       
       
\usepackage{enumerate}       
\usepackage{mathabx}
\usepackage[active]{srcltx}
\usepackage{paralist} 
\numberwithin{equation}{section}       
       
 \makeatletter      
  \@namedef{subjclassname@2020}{%
  \textup{2020} Mathematics Subject Classification}    
   \makeatother    
       
\usepackage{t1enc}       
\usepackage[utf8x]{inputenc}

\newtheorem{theorem}{Theorem}[section]        
\newtheorem{lemma}[theorem]{Lemma}

\newtheorem{proposition}[theorem]{Proposition}

\theoremstyle{definition}       
\newtheorem{definition}[theorem]{Definition}       
\theoremstyle{remark}       
         
\newcommand{\mc}[1]{\mathcal{#1}}       
\newcommand{\ms}[1]{\mathscr{#1}}       
\newcommand{\mbb}[1]{\mathbb{#1}}

\newcommand{\setm}{\setminus}       
\newcommand{\empt}{\emptyset}       
\newcommand{\subs}{\subset}       
\newcommand{\dom}{\operatorname{dom}}       
\newcommand{\ran}{\operatorname{ran}}       
\def\<{\left\langle}       
\def\>{\right\rangle}

\newcommand{\clsco}{ClScat}
\newcommand{\dio}{Disc}
\newcommand{\sco}{Scat}
\newcommand{\nwdo}{Nwd}

\newcommand{\dioo}{Disc_{\omega}}
\newcommand{\scoo}{Scat_{\omega}}
\newcommand{\nwdoo}{Nwd_{\omega}}

\newcommand{\tocn}[2]{#1^{(#2)}}

\author[I. Juh\'asz]{Istv\'an Juh\'asz}       
\address       
      { Alfr{\'e}d R{\'e}nyi Institute of Mathematics, HUN-REN}
       
\email{juhasz@renyi.hu}

\author[L. Soukup]{Lajos Soukup}       
\address       
      { Alfr{\'e}d R{\'e}nyi Institute of Mathematics,  HUN-REN
}       
\email{soukup@renyi.hu}       
       
\author[Z. Szentmikl\'ossy]{Zolt\'an Szentmikl\'ossy}       
\address{Eötvös University of Budapest}       
\email{szentmiklossyz@gmail.com}

\title{Selectively pseudocompact spaces}       
\thanks{The preparation of this paper was
supported by  OTKA grant   K129211}       
\date{\today}       
       

\begin{document}

\begin{abstract}       
    A novel {\em selection principle} was introduced  by Dorantes-Aldama and Shakhmatov:
    a topological space $X$ is termed {\em selectively pseudocompact} if for any 
    sequence $\<U_n:n\in {\omega}\>$ of pairwise disjoint non-empty open sets of $X$, one can choose points $x_n\in U_n$ such that 
    the sequence $\<x_n:n\in {\omega}\>$ has an accumulation point. 
    
 In this paper, we explore various versions of this principle when we permit the 
 selection of finite, scattered,  or nowhere dense sets instead of just singletons.

 We develop a method to prove that the  aforementioned  versions of
 selective pseudocompactness  
are indeed distinct from one another. 
\end{abstract}       
       
\maketitle

\section{Introduction}

A novel {\em selection principle} was introduced in \cite[Definition 2.2]{DoSh17}:
a topological space $X$ is termed {\em selectively pseudocompact} if for any 
sequence $\<U_n:n\in {\omega}\>$ of pairwise disjoint non-empty open sets of $X$, one can choose points $x_n\in U_n$ such that 
the sequence $\<x_n:n\in {\omega}\>$ has an accumulation point.

A selectively pseudocompact space is evidently pseudocompact. However, the converse is not true 
since  there exist pseudocompact spaces where all their countable subsets are closed 
(see \cite[Exaple 2.7]{DoSh17}  or  \cite[Theorem 2] {Sa86}).

In \cite{DoSh17} and \cite{DoSh20} the authors introduced and  investigated several variations 
of the selective pseudocompactness property. Yet, it transpires that all these variations are 
more stringent than selectively pseudocompactness itself.

When considering a selection principle, a natural query arises: 
What transpires  if we alter  the requirements concerning the chosen sets? 
In the original definition of selective pseudocompactness, one had    to select singletons. 
What unfold if we permit the selection of  finite sets rather than singletons? Does this engender 
a distinct class of spaces? This query marked  the  inception  of our research.

\begin{definition}\label{df:f-sel-pseudo}
    Let $X$ be a topological space and $\mc A\subs \ms P(X)$.
    We say that  $X$ is {\em $\mc A$-selectively pseudocompact} ($\mc A$-SP,  in short)  if given any 
    sequence $\<U_n:n\in {\omega}\>$ of  non-empty open sets, one can pick sets  
    $A_n\in \mc A\cap \mc P(U_n)$ such that 
    the family  $\{A_n:n\in {\omega}\}$ has an accumulation point (i.e., it is not locally finite). 
    \end{definition}
    Hence, $X$ is selectively pseudocompact iff it is ${[X]}^{1}$-SP.

If every non-empty open subset of $X$ contains an element of $\mc A$, then 
 the $\mc A$-selectively pseudocompact property is   between pseudocompactness and 
selectively pseudocompactness.
Since a scattered space is pseudocompact iff 
the set of isolated points is relatively countable compact,
from now on we will consider only dense-in-itself spaces.  

The most natural variation of selectively pseudocompactness is when we allow selecting finite
sets instead of singletons.
 Let us say that a space $X$ is {\em fin-selectively pseudocompact}   ({\em fin-SP}, in short) 
 iff $X$ is 
${[X]}^{<{\omega}}$-SP.  
In other versions one can allow selecting discrete, scattered or nowhere dense subsets as well.   
To formulate the corresponding notions
let us fix first the following terminology.
\begin{definition}\label{df:ope}
If $X$  is a topological space, write
\begin{align*}
    \dio(X)=&\{Y\subs X:\text{$Y$ is discrete}\},\\
    \clsco(X)=&\{Y\subs X:\text{$Y$ is closed and scattered}\},\\
\sco(X)=&\{Y\subs X:\text{$Y$ is scattered}\},\\
\nwdo(X)=&\{Y\subs X:\text{$Y$ is nowhere dense}\}.
\end{align*}
To shorten the formulation of the results, we 
say that $X$ is {\em $\dio$-SP} iff $X$ is $\dio(X)$-SP. 
The notions $\clsco$-SP,  $\sco$-SP and $\nwdo$-SP are introduced similarly. 
\end{definition}

If $X$ is a topological space and $\mc A\subs \mc P(X)$, write 
$\overline{\mc A}=\{\overline A:A\in \mc A\}.
$
If $X$ is crowded, then $$\clsco(X)\subseteq \overline{\sco(X)}=\overline{\dio(X)}
\subseteq \overline{\nwdo(X)},$$ but  
$$\clsco(2^{\omega})\ne \overline{\sco(2^{\omega})}.$$

The following lemma is straightforward.  
\begin{lemma}\label{lm:closure}
    Assume that $X$ is a regular space and $\mc A\subs \mc P(X)$.
    Then $X$ is $\mc A$-SP iff it is $\overline{\mc A}$-SP.
  \end{lemma}

  \begin{proposition}\label{cr:closure}
    A regular space $X$ is $\sco$-SP iff it is $\dio$-SP.
    \end{proposition}
        \begin{proof}
        Apply Lemma \ref{lm:closure} twice:  
            $X$ is  $\sco(X)$-SP iff $X$ is   $\overline{\sco(X)}$-SP, and 
        $X$ is   $\dio(X)$-SP iff $X$ is  $\overline{\dio(X)}$-SP.
        Finally, observe that  $\overline{\dio(X)}=\overline{\sco(X)}$.         
        \end{proof}

Let $X$ be a crowded space and consider the following statements:
\begin{inparaenum}[(i)]
\item $X$ is $1$-SP,
\item $X$ is fin-SP,
\item $X$ is $\clsco$-SP,
\item $X$ is $\sco$-SP,
\item $X$ is $\nwdo$-SP,
\item $X$ is pseudocompact. 
\end{inparaenum} 

Clearly, (i)$\implies$ (ii)$\implies$(iii)$\implies$(iv)$\implies$ (v)$\implies$ (vi). 
In the next section we prove the   Theorem \ref{tm:main-3} below
which shows that the implications, apart from the implication (iv)$\implies$ (v) are not reversible.

\begin{theorem}[CH]\label{tm:main-3}
    There are countably tight, 0-dimensional Hausdorff spaces $X_i$ for $i=1,2,3,4$ such that 
    \begin{enumerate}[(1)]
    \item $X_1$ is $fin$-SP, but not $1$-SP,
    \item $X_2$ is $\clsco$-SP, but not $fin$-SP,
    \item $X_3$ is $\sco$-SP, but not $\clsco$-SP, 
    \item $X_4$ is pseudocompact, but not $\nwdo$-SP. 
    \end{enumerate}
    \end{theorem}

The method we apply in the proof of this  Theorem is not suitable 
to separate properties (iv) and (v)
because the constructed spaces contain dense, crowded open metric subspaces,
so by Proposition \ref{cr:dense-c} below, the constructed spaces can not separate these properties.

\begin{proposition}\label{cr:dense-c}
    Assume that a crowded space $X$ contains a metric, open dense subset $D$. 
    Then $X$ is $Discrete(X)$-SP  iff it is $Nwd(X)$-SP.
\end{proposition}

\begin{proof}
    Let $\{U_n:n<{\omega}\}$  be a sequence of pairwise disjoint non-empty open sets. 
    Since $D$ is dense and open, we can assume that $U_n\subs D$.
 
Since $X$ is $\nwdo$-SP
we can pick  $A_n\in \nwdo(X)\cap \mc P(U_n)$ such that $\{A_n:n<{\omega}\}$ has an accumulation point $z$. 

Since $D$ is metric, for each $n\in {\omega}$ there is a discrete $B_n\subs D$ with $A_n\subs \overline{B_n}$.
We can clearly assume that $B_n\subs U_n$. 
Then $z$ is an accumulation point of  $\{B_n:n<{\omega}\}$.

Thus $X$ is $\dio(X)$-SP.
%
%
%
\end{proof}

Finally ,  in Theorem \ref{tm:sc-nwd}, using a different method, we can prove that a version of  
implication (iv)$\implies$ (v) is not reversible, as well.  
Given a topological space $X$, write
\begin{align*}
\scoo(X)=&\{B\in {[X]}^{{\omega}}:\text{$B$ is scattered}\},\\
\nwdoo(X)=&\{S\in {[X]}^{{\omega}}:\text{$S$ is nowhere dense}\}.
\end{align*}

\begin{theorem}[CH]\label{tm:sc-nwd}
    There is a crowded, 0-dimensional Hausdorff topological space $X$ which is  
    $\nwdoo(X)$-SP, but not $\scoo(X)$-SP.
    \end{theorem}


\section{General separation result}

\subsection*{Preliminary definitions}

\begin{definition}\label{df:base}
A {\em base-enumeration $\mc B$} is a sequence $\<B_i:i<{\xi}\>$ for some ordinal
such that $\{B_i:i<{\xi}\}$ is a clopen base of a (not necessarily $\mbb T_2$) topology 
${\tau}$ on some set $X$.

If $\mc B^{\nu}$ is a base-enumeration, we use the notation ${\xi}^{\nu}$,
$B^{\nu}_i$ for $i<{\xi}^{\nu}$, ${\tau}^{\nu}$ and $X^{\nu}$.

If $\mc B^0$ and $\mc B^1$ are base-enumerations, write $\mc B^0\preceq \mc B^1$
iff 
\begin{enumerate}[(1)]
\item $X^0\subs X^1$, ${\xi}^0\le {\xi}^1$, 
\item $B^0_i=B^1_i\cap X_0$ for $i<{\xi}^0$,
\item $B^0_i\subs B^0_j$ iff $B^1_i\subs B^1_j$ for $i,j<{\xi}^0$ ,
\item $B^0_i\cap B^0_j=\empt$ iff $B^1_i\cap B^1_j=\empt$ for $i,j<{\xi}^0$. 
\end{enumerate}
\end{definition}

\begin{lemma}\label{lm:base-enumeration-limit}
Assume that  ${\mu}$ is a limit ordinal, and $\<\mc B^{\nu}:{\nu}<{\mu}\>$ is a 
$\preceq$-increasing sequence of base-enumerations.
Write 
\begin{enumerate}[(i)]
\item ${\xi}^{\mu}=\sup_{{\nu}<{\mu}}{\xi}^{\nu}$,
\item $B^{\mu}_i=\bigcup\{B^{\nu}_i:{\nu}<{\mu}, i<{\xi}^{\nu}\}$ for $i<{\xi}^{\mu}$,
\end{enumerate}
Then $\mc B^{\mu}=\<B^{\mu}_i:i<{\xi}^{\mu}\>$ is a base-enumeration, 
and $\mc B^{\nu}\preceq \mc B^{\mu}$ for each ${\nu}<{\mu}$.

Moreover, if ${\tau}^{\nu}$ is $\mbb T_2$ for cofinally many ${\nu}<{\mu}$,
then ${\tau}_{\mu}$ is also ${\tau}_2$.

Finally, if $a\in X^{\sigma}$ and  $A\subs X^{\sigma}$ for some 
${\sigma}<{\mu}$ such that $a\in \overline{A}^{{\tau}_{\nu}}$ for each 
${\sigma}<{\nu}<{\mu}$, then $a\in \overline{A}^{{\tau}_{\mu}}$
\end{lemma}
We will write $\lim_{{\nu}\to {\mu}}\mc B^{\nu}$ for $\mc B^{\mu}$.

\begin{proof}
Straightforward.
\end{proof}

\subsection*{The separation result}

For $n\in {\omega}$ write $C_n=\{n\}\times 2^{\omega}$, and let $C=\bigcup_{n<{\omega}}C_n$.
Denote ${\varepsilon}$ the natural, Euclidean topology on $C$. 

Given any functions  ${\varphi}:2^{\omega}\to Z$ 
and ${\psi}:\mc P(2^{\omega})\to Z$, 
and a natural number $n$, define the functions $\tocn{{\varphi}}{n}:C_n\to Z$
and $\tocn{{\psi}}{n}:\mc P(C_n)\to Z$
by the formulas $\tocn{{\varphi}}{n}(\<n,x\>)={\varphi}(x)$
and $\tocn{{\psi}}{n}(A)={\psi}(\{x:\<x,n\>\in A\})$.

Let 
\begin{displaymath}
\mbb E=\{e:\text{$e$ is a function, $\dom(e)\in {[{\omega}]}^{{\omega}}$, and $\empt\ne e(n)\subs C_n$ for each $n\in \dom(e)$}\},
\end{displaymath}
and 
\begin{displaymath}
\mbb O=\{o\in \mbb E:\text{$o(n)$ is a clopen subset of $C_n$ for each $n\in \dom(o)$}\}.
\end{displaymath}
For $e\in \mbb E$ let $$R(e)=\bigcup \ran(e).$$
For $f,g\in \mbb E$ write 
\begin{displaymath}
    \text{$f\sqsubset g$  iff $\dom(f)\subs \dom(g)$ and $f(n)\subs g(n)$ for each $n\in \dom(f)$},
\end{displaymath}
\begin{displaymath}
    \text{$f\sqsubset^* g$  iff $\dom(f)\subs^* \dom(g)$ and $|\{n\in\dom(f):f(n)\not\subs g(n)\}|<{\omega}$},
\end{displaymath}
and 
\begin{displaymath}
    \text{$f\bot g$  iff $|\{n\in\dom(f)\cap \dom(g):f(n)\cap g(n)\ne \empt\}|<{\omega}$}.
\end{displaymath}
For $h,o\in \mc E$ define the function $h\dotdiv o $ as follows:
 $$\dom(h\dotdiv o)=\{n\in \dom(h): h(n)\not\subs o(n)\},$$
 and 
 \begin{displaymath}
 {(h\dotdiv o)}(n)=\left\{\begin{array}{ll}
 {h(n)\setm o(n)}&\text{if $n\in \dom(o)$,}\\
 {h(n)}&\text{if $n\notin \dom(o)$}.
 \end{array}\right.
 \end{displaymath}

\begin{definition}\label{df:gh}
If ${\mu}:\mc P(2^{\omega})\to {\omega}\cup\{\infty\}$ and $M\in {\omega}$, define 
the families $\mbb G({\mu},M)$ and $\mbb H({\mu})$ and follows
\begin{displaymath}
\mbb G({\mu},M)=\{g\in \mbb E: \tocn{\mu}n(g(n))\le M \text{ for each $n\in \dom(g)$}\}
\end{displaymath}
and 
\begin{displaymath}
\mbb H({\mu})=\{h\in \mbb E: \lim_{n\in \dom(h)}\tocn{\mu}n(h(n))=\infty\}.
\end{displaymath}
\end{definition}

\begin{definition}\label{df:shrink-dense-large}
Let $\mbb G,\mbb H\subs \mbb E$.
\begin{enumerate}[(1)]
\item  $\mbb H$ is {\em dense} iff for each $o\in \mbb O$ there is $h\in \mbb H$ with $h\sqsubset o$.
\item $\mbb H$ is {\em shrinkable } iff for each $h\in \mbb H$ 
and $\mc O\in {[\mbb O]}^{{\omega}}$
there is $h'\in \mbb H$  such that
\begin{displaymath}
    h'\sqsubset h \text{ and }\forall o\in \mc O \ \big( h'\sqsubseteq^* o\ \lor \ h'\bot o\big). 
\end{displaymath}
\item $\mbb H$ is {\em $\mbb G$-large} iff for each $g\in \mbb G$ and $\mc H\in {[\mbb H]}^{{\omega}}$
there is $o\in \mbb O$ such that 
\begin{displaymath}
g\sqsubset o \text{ and } h\dotdiv o\in \mbb H \text{ for each }h\in \mc H. 
\end{displaymath}
\end{enumerate}
\end{definition}

\begin{theorem}\label{tm:main-1}
Assume that $M\in {\omega}$ and ${\mu}:\mc P(2^{\omega})\to {\omega}\cup \{\infty\}$ is a monotone, subadditive function 
such that \begin{enumerate}[(i)]
\item ${\mu}(U)=\infty$ for each non-empty clopen $U\subs 2^{\omega}$, 
\item for each $a,b\subs 2^{\omega}$, if ${\mu}(a)\le M<{\mu}(b)$, then 
\begin{displaymath}
\sup\{{\mu}(b\setm U)+M:a\subs U\in clopen(2^{\omega})  \}\ge {\mu}(b).
\end{displaymath}
\end{enumerate}
Then $\mbb H({\mu})$ is dense, shrinkable and $\mbb G({\mu},M)$-large.
\end{theorem}

\begin{theorem}[CH]\label{tm:main-2}
    Assume that $\mbb G,\mbb H \subs \mbb E$ such that 
    $\mbb H$ is dense, shrinkable and $\mbb G$-large.
    Then, there is a 0-dimensional $\mbb T_2$ topological space $X=\<C\cup {\omega}_1,{\tau}\>$
    such that 
    \begin{enumerate}[(i)]
    \item $\mbb T(X)={\omega}$  and ${\omega}_1$ is left-separated in its natural order, 
    \item $C$ is a dense open subspace of $X$ and   ${\tau}\restriction C={\varepsilon}$,
    \item $\overline{R(g)}\cap {\omega}_1=\empt$ for each $g\in \mbb G$,
    \item $\overline{R(h)}\cap {\omega}_1\ne \empt$ for each $h\in \mbb H$.
    \end{enumerate}
    \end{theorem}

    \begin{proof}[Proof of Theorem \ref{tm:main-1}]
Since ${\mu}(U)=\infty$ for each non-empty clopen subset of $C$, 
we have $\mbb O\subs \mbb H({\mu})$, so $\mbb H({\mu})$ is dense.

To show that $\mbb H({\mu})$ is shrinkable, fix 
$h\in \mbb H$ and $\{o_m:m<{\omega}\}\in {[\mbb O]}^{{\omega}}$.

Define $h_m\in \mbb H({\mu})$
and ${\chi}_m\in 2$ for $m<{\omega}$ such that 
\begin{displaymath}
h=h_0\sqsupset h_1\sqsupset\dots
\end{displaymath}
and writing $I_m=\dom(h_m)$ we have 
\begin{enumerate}[(i)]
\item $\tocn{\mu}m(h_m(i))\ge m$ for each  $i\in I_m$,
\item if ${\chi}_m=1$ then $h_{m+1}(i)=h_m(i)\cap o_m(i)$ for each $i\in I_{m+1}$,

\item if ${\chi}_m=0$ then  
either $\dom(o_m)\cap \dom(h_{m+1})=\empt$ 
and $h_{m+1}(i)=h_m(i)$ for each $i\in I_{m+1}$, 
or
$h_{m+1}(i)=h_m(i)\setm  o_m(i)$ for each $i\in I_{m+1}$. 
\end{enumerate}

Assume that $h_m$ and so $I_m$ are constructed. 

If $\dom(h_m)\setm \dom(o_m)$ is infinite, then let 
${\chi}_m=0$ and 
then pick $K_{m}\in {[\dom(h_m)\setm \dom(o_m)]}^{{\omega}}$ such that 
${\mu}(h_m(i))\ge m+1$ for each $i\in K_m$ and let 
$h_{m+1}=h_m\restriction K_m$.

If $\dom(h_m)\setm \dom(o_m)$ is not infinite,
then $J_m=\dom(h_m)\cap \dom(o_m)$ is infinite. 
Since $\lim_{i\in J_m}\tocn{\mu}m(h_m(i))=\infty$ and ${\mu}$ is subadditive,  
either 
\begin{enumerate}[(a)]
\item $\limsup_{i\in J_m}\tocn{\mu}m(h_m(i)\cap o_m(i))=\infty$,
\end{enumerate}
or 
\begin{enumerate}[(b)]
    \item $\limsup_{i\in J_m}\tocn{\mu}m(h_m(i)\setm o_m(i))=\infty$
    \end{enumerate}
    holds. 

    If (a) holds, then we can find $K_n\in {[J_n]}^{{\omega}}$
    such that $\tocn{\mu}m(h_m(i)\cap o_m(i))\ge m+1$ for each $i\in K_m$
and $\lim_{i\in K_m}\tocn{\mu}m(h_m(i)\cap o_m(i))=\infty$.
Then let ${\chi}_m=1$, $\dom(h_{m+1})=K_m$ and 
$h_{m+1}(i)=h_m(i)\cap o_m(i)$.

If (b) holds, then we can find $K_n\in {[J_n]}^{{\omega}}$
such that $\tocn{\mu}m(h_m(i)\setm  o_m(i))\ge m+1$ for each $i\in K_m$
and $\lim_{i\in K_m}\tocn{\mu}m(h_m(i)\setm o_m(i))=\infty$.
Then let ${\chi}_m=0$, $\dom(h_{m+1})=K_m$ and 
$h_{m+1}(i)=h_m(i)\setm o_m(i)$.

Finally, let 
\begin{displaymath}
a_m=\min(\dom(h_m)\setm\{a_\ell:\ell<m\}.)
\end{displaymath}
for $m\in {\omega}$, and define the function $h'$ as follows: 
\begin{enumerate}[$\bullet$]
\item $\dom(h')=\{a_m:m<{\omega}\}$,
\item $h'(a_m)=h_m(a_m)$,
\end{enumerate}
Since $\tocn{\mu}m(h_m(a_m))\ge m$, we have $h'\in \mbb H({\mu})$.

By the construction, 
\begin{enumerate}[$\bullet$]
\item if ${\chi}_m=1$, then $h'(a_i)\subs o_m(a_i)$ for $i>m$, and 
\item if ${\chi}_m=0$, then $h'(a_i)\cap o_m(a_i)=\empt$ for $i>m$.  
\end{enumerate}
So either $h'\sqsubset^* o_m$ or $h'\bot o_m$.
Thus, $\mbb H({\mu})$ is shrinkable. 

\medskip

To show that $\mbb H({\mu})$  is $\mbb G({\mu},M)$-large fix  $g\in \mbb G$ and 
$\mc H=\{h_m:m\in {\omega}\}\in {[\mbb H]}^{{\omega}}$.

Fix a sequence $\{n_k\}_{k\in {\omega}}$ of natural numbers which converges to infinity such that 
\begin{displaymath}
\forall k\in {\omega}\ \forall i \ge n_k\  \forall \ell\le k \
\text{if $i\in \dom(h_\ell)$, then }
\tocn{\mu}\ell(h_\ell(i))\ge k. 
\end{displaymath}

If $n_k\le i<n_{k+1}$ for each $\ell\le k$,
if $i\in \dom(g)\cap\dom(h)\ell$, then 
by (ii) we can pick clopen $U^\ell_i\subs C_i$ such that
$g(i)\subs U^\ell_i$ and $\tocn{\mu}\ell(h_\ell(i)\setm U^\ell_i)\ge k$. 

Define $o\in \mbb O$ such that
$\dom(o)=\dom(g)$ and  
 $$o(i)=\bigcap\{U^\ell_i:{\ell\le k}, i\in \dom(h_\ell)\}.$$ 
Then $g(i)\subs o(i)$ and ${\mu}^{(\ell)}((h_\ell\dotdiv o)(i))\ge k$.
Thus, $\lim_{i\in \dom(h_\ell)}\tocn{\mu}\ell(h_\ell(i)\setm o(i))=\infty$, and so 
$h_\ell\dotdiv o\in \mbb H({\mu})$.
\end{proof}

\begin{proof}[Proof of Theorem \ref{tm:main-2}]
We will construct our space by transfinite recursion. To do so, we need the following notion. 
    
We say that a triple $\mbb T=\<C\cup {\gamma},\mc B, \mc F\>$
is a {\em good triple} iff
\begin{enumerate}[(a)]
    \item ${\gamma}\le {\omega}_1$ is an  ordinal,
    \item    $\mc B=\{B_i:i<{\xi}\}$ for some ordinal ${\xi}\le {\omega}_1$,
 \item $\mc B$ is a base of a 0-dimensional $\mbb T_2$ topology 
${\tau}$ on  $X=C\cup {\gamma}$,
\item ${\tau}\restriction C={\varepsilon}$, and  
$\{B_i:i<{\omega}\}$ is a base of  ${\varepsilon}$,
\item for each ${\omega}\le {\zeta}<{\xi}$ there is $o_{\zeta}\in \mbb O$
such that $B_{\zeta}\cap C=R(o_{\zeta})$,
\item $C$ is dense open in ${\tau}$,
\item $\mc F\in {[{\gamma}\times \mbb H]}^{\le {\omega}}$,
\item  If $\<{\alpha},h\>\in \mc F $ 
then ${\alpha}\in \overline{R(h)}^{\tau}$.   
\end{enumerate}
We say that $\mbb T$ is {\em countable} iff both ${\gamma}$ and ${\xi}$ are countable ordinals.

If $\mbb T_{\nu}$ is a good triple, we use the notation 
${\gamma}^{\nu}$, ${\tau}^{\nu}$,  $X^{\nu}$, ${\xi}^{\nu}$, $B^{\nu}_i$ for 
$i<{\xi}^{\nu}$ and $h^{\nu}_j$ for $j<{\gamma}^{\xi}$.

If $\mbb T^0$ and $\mbb T^1$ are good triples, write $\mbb T^0\preceq T^1$ iff
\begin{enumerate}[(A)]
\item ${\gamma}^0\le {\gamma}^1$, ${\xi}^0\le {\xi}^1$,
\item $B^0_i=B^1_i\cap X^0$ for $i<{\xi}^0$,
\item $B^0_i\subs B^0_k$ iff  
 $B^1_i\subs B^1_k$ for $i,k<{\xi}^0$,  
\item $B^0_i\cap B^0_k=\empt$ iff  
 $B^1_i\cap B^1_k=\empt$ for $i,k<{\xi}^0$,
 \item ${\tau}^0\restriction {\gamma}_0={\tau}^1\restriction {\gamma}_0$   
\item $\mc F^0\subs \mc F^1$.
\end{enumerate}

\begin{lemma}\label{lm:1}
If $\mbb T^0$ is a countable good triple and   $h\in \mc H$, then 
there is a countable good triple $\mbb T^1$ such that $\mbb T^0\preceq T^1$ and  
$\<j,h\>\in \mc F^1$ for some $j<{\gamma}^1$.  
\end{lemma}

\begin{proof}
Fix $o_{i}\in \mbb O$ with $B^0_{i}\cap C=R(o_{i})$ for $i<{\xi}^0$.
Let $\mc O=\{o_i:i<{\xi}^0\}$.
Since $\mbb H$ is shrinkable, there is $h'\in \mbb H$  such that
\begin{displaymath}
    h'\sqsubset h \text{ and }\forall o\in \mc O \ \big( h'\sqsubseteq^* o_i\ \lor \ h'\bot o_i\big). 
\end{displaymath}
for each $i<{\xi}^0$. 
Let ${\gamma}^1={\gamma}^0+1$ and 
${\xi}^1={\xi}^0+{\omega}$. For $i<{\xi}^0$ let 
\begin{displaymath}
{B^1_i}=\left\{\begin{array}{ll}
{B^0_i}&\text{if $h'\bot o_i$,}\\\\
{B^0_i\cup\{{\gamma}^0\}}&\text{if $h'\sqsubseteq^* o_i$,}\\
\end{array}\right.
\end{displaymath}
For $\ell<{\omega}$ let
\begin{displaymath}
B^1_{{\xi}^0+\ell}=\{{\gamma}^0\}\cup \bigcup\{h'(m):m\in \dom(h')\setm  \ell\}.
\end{displaymath}
Write $\mc B^1=\{B^1_i:i<{\xi}^1\}$.
Then the good tripe $\mbb T^1=\<C\cup {\xi}^1,\mc B^1,\mc F^0\cup\{\<h,{\gamma}_0\>\}\>$
meets the requirements. 
\end{proof}

\begin{lemma}\label{lm:2}
If $\mbb T^0$ is a countable good triple and   $g\in \mbb G$, then 
there is a countable good triple $\mbb T^1$ such that  $\mbb T^0\preceq T^1$ and 
$\overline{R(g)}^{{\tau}^1}\cap {\gamma}^0=\empt$. 
\end{lemma}

\begin{proof}
Fix $o_{i}\in \mbb O$ with $B^0_{i}\cap C=R(o_{i})$ for $i<{\xi}^0$.
Write $\mc O=\{o_i:i<{\xi}^0\}$.
Let $\mc H=\{h:\<{\alpha},h\>\in \mc F^0\}$.
Since  $\mbb H$ is {\em $\mbb G$-large},  $\mc H\in {[\mbb H]}^{{\omega}}$
there is $o\in \mbb O$ such that 
\begin{displaymath}
g\sqsubset o \text{ and } h\dotdiv o\in \mbb H \text{ for each }h\in \mc H\cup \mc O. 
\end{displaymath}
Let ${\gamma}^1={\gamma}^0$. If ${\xi}={\zeta}\dotplus {\omega}$, then  
${\xi}^1={\xi}\dotplus {\zeta}$, $B^1_i=B^0_i$ for $i<{\xi}$ and 
$B^1_{{\xi}\dotplus {\eta}}=B^0_{{\omega}\dotplus {\eta}}\setm R(o)$.
Write $\mc B^1=\{B^1_i:i<{\xi}^1\}$.

Then the good tripe $\mbb T^1=\<C\cup {\gamma}^1,\mc B^1,\mc F^0\>$
meets the requirements. 
\end{proof}

After this preparation we can carry out the following construction. 
By transfinite induction on ${\nu}\le {\omega}_1$
define a $\preceq$-increasing sequence of goof triple, $\<T^{\nu}:{\nu}\le {\omega}_1\>$
as follows.

Let $\{f_{\xi}:{\xi}<{\omega}_1\}$ be an ${\omega}_1$-abundant enumeration of 
$\mbb G\cup \mbb H$. 

Let ${\gamma}^0=0$, let $\<B^0_i:i<{\omega}\>$ be a clopen base of $C$ in the Euclidean topology,
and let $\mc F^0=\empt$. Then $\mbb T^0=\<C,\mc B^0,\empt\>$ is a good countable triple. 

If ${\nu}$ is a limit ordinal, and $\mbb T^{\eta}$ is defined for ${\eta}<{\nu}$, take the natural ``union'', 
i.e. let  ${\gamma}^{\nu}=\sup\{{\gamma}^{\eta}:{\eta}<{\nu}\}$,
 ${\xi}^{\nu}=\sup\{{\xi}^{\eta}:{\eta}<{\nu}\}$,
 $\mc F^{\nu}=\bigcup \{{\mc F}^{\eta}:{\eta}<{\nu}\}$, $\mc B^{\nu}=\{B^{\nu}_i:i<{\xi}^{\nu}\}$,
 where
\begin{displaymath}
B^{\nu}_i=\bigcup\{B^{\eta}_i:i<{\xi}^{\eta}\}.
\end{displaymath}
Then  $\mbb T^{\nu}=\<C\cup {\gamma}^{\nu},\mc B^{\nu},\mc F^{\nu}\>$ 
is a good triple, and $\mc T^{\eta}\prec \mc T^{\nu}$ for ${\eta}<{\nu}$.

Finally, assume that ${\nu}={\mu}+1$.

If $f_{{\mu}}\in \mbb G$, then applying Lemma \ref{lm:2}
for $\mbb T^{\mu}$ and $f_{\xi}$ we obtain $\mbb T^{\nu}$ such that 
$\overline{R(f_{\mu})}^{{\tau}^{\nu}}\cap {\gamma}^{\mu}=\empt$. 

If $f_{{\mu}}\in \mbb H$, then applying Lemma \ref{lm:1}
for $\mbb T^{\mu}$ and $f_{\xi}$ we obtain $\mbb T^{\nu}$ such that 
$\<{\xi},f_{\mu}\>\in \mc F^{\nu}$.  

We claim that $X=\<C\cup{\omega}_1,{\tau}^{{\omega}_1}\>$ satisfies the requirements. 

The construction clearly guarantees that
\begin{enumerate}[(i)]
    \addtocounter{enumi}1    
    \item $C$ is a dense open subspace of $X$ and   ${\tau}\restriction C={\varepsilon}$,
     
\item $\overline{R(g)}\cap {\omega}_1=\empt$ for each $g\in \mbb G$,
\item $\overline{R(h)}\cap {\omega}_1\ne \empt$ for each $h\in \mbb H$.
\end{enumerate} 
It is also clear that ${\omega}_1$ is left separated in the natural order. 

TO  check that $t(X)={\omega}_1$ it is enough 
to show that $t(\<{\omega}_1,{\tau}^{{\omega}_1}\>)={\omega}$.

Assume that ${\alpha}\in {\omega}_1,$ $A\in {[{\omega}_1]}^{{\omega}_1}$
and ${\alpha}\in \overline{A}^{{\tau}^{{\omega}_1}}$.

Using a simple closing argument we can find ${\nu}<{\omega}_1$
 such that ${\alpha}<{\gamma}^{\nu}$  and 
$
 {\alpha}\in \overline{A\cap {\gamma}_{\nu}}^{{\tau}_{\nu}}
$
. Then ${\alpha}\in \overline{A\cap {\gamma}_{\nu}}^{{\tau}_{{\omega}_1}}$
because 
$
    {\tau}^{{\omega}_1}\restriction {\gamma}^{\nu}={\tau}^{\nu}\restriction {\gamma}^{\nu}
$
by (E).
\end{proof}

\subsection*{Applications}
\begin{proof}[Proof of Theorem \ref{tm:main-3}]
(1) Let $M=0$ and define ${\mu}_1:\mc P(C)\to {\omega}\cup\{\infty\}$ by the formula 
\begin{displaymath}
{{\mu}_1}(A)=\left\{\begin{array}{ll}
{|A|}&\text{ if $A$ is finite, }\\
{\infty}&\text{ if $A$ is infinite}.
\end{array}\right.
\end{displaymath}

Applying Theorem \ref{tm:main-1}
we obtain that $\mbb H({\mu})$ is dense, shrinkable and $\mbb G({\mu},M)$-large.
So we can apply Theorem \ref{tm:main-2} for 
$\mbb H({\mu})$ and $\mbb G({\mu},M)$ to obtain a space $X_1$
which satisfies \ref{tm:main-2}.(i)-(iv). 
The sequence $\<C_n:n<{\omega}\>$ witnesses that 
$X_1$ is not ${1}$-SP. 

\medskip
\noindent
(2) Let $M=1$ and define ${\mu}_2:\mc P(C)\to {\omega}\cup\{\infty\}$ by the formula 
\begin{displaymath}
{{\mu}_2}({A})=\left\{\begin{array}{ll}
{ht(\overline{A})}&\text{if $\overline{A}$ is scattered of finite height,}\\
{\infty}&\text{otherwise}.
\end{array}\right.
\end{displaymath}
Applying Theorem \ref{tm:main-1}
we obtain that $\mbb H({\mu})$ is dense, shrinkable and $\mbb G({\mu},M)$-large.
So we can apply Theorem \ref{tm:main-2} for 
$\mbb H({\mu})$ and $\mbb G({\mu},M)$ to obtain a space $X_2$
which satisfies \ref{tm:main-2}.(i)-(iv). 
The sequence $\<C_n:n<{\omega}\>$ witnesses that 
$X_2$ is not ${fin}$-SP.

\medskip
\noindent
(3) Let $M=1$ and define ${\mu}_3:\mc P(C)\to {\omega}\cup\{\infty\}$ by the formula 
\begin{displaymath}
{{\mu}_3}(A)=\left\{\begin{array}{ll}
{0}&\text{if $\overline{A}$ is scattered,}\\
{\infty}&\text{otherwise}.
\end{array}\right.
\end{displaymath}
Applying Theorem \ref{tm:main-1}
we obtain that $\mbb H({\mu})$ is dense, shrinkable and $\mbb G({\mu},M)$-large.
So we can apply Theorem \ref{tm:main-2} for 
$\mbb H({\mu})$ and $\mbb G({\mu},M)$ to obtain a space $X_3$
which satisfies \ref{tm:main-2}.(i)-(iv). 
The sequence $\<C_n:n<{\omega}\>$ witnesses that 
$X_3$ is not $\clsco$-SP.

\medskip
\noindent
(4) Let $M=1$ and define ${\mu}_4:\mc P(C)\to {\omega}\cup\{\infty\}$ by the formula 
\begin{displaymath}
{{\mu}_4}(A)=\left\{\begin{array}{ll}
{0}&\text{if $A$ is nowhere dense,}\\
{\infty}&\text{if $int(\overline A)\ne \empt$}.
\end{array}\right.
\end{displaymath}
Applying Theorem \ref{tm:main-1}
we obtain that $\mbb H({\mu})$ is dense, shrinkable and $\mbb G({\mu},M)$-large.
So we can apply Theorem \ref{tm:main-2} for 
$\mbb H({\mu})$ and $\mbb G({\mu},M)$ to obtain a space $X_4$
which satisfies \ref{tm:main-2}.(i)-(iv). 
The sequence $\<C_n:n<{\omega}\>$ witnesses that 
$X_4$ is not ${\nwdo}$-SP.

\end{proof}

\section{Scattered-SP vs nwd-SP}


Instead of Theorem \ref{tm:sc-nwd}
we prove the following stronger result. 
\begin{theorem}[CH]\label{tm:sc-nwdplus}
There is a crowded, 0-dimensional Hausdorff topological space $X$ such that $|X|=w(X)={\omega}_1$,
every $A\in \dioo(X)$ is closed, and $X$ is $\nwdoo(X)$-SP.
\end{theorem}

\begin{proof}[Proof of Theorem \ref{tm:sc-nwdplus}]
    We say that a triple $\mbb T=\<{\gamma},\mc B, \mc F\>$
    is a {\em nice triple} iff
    \begin{enumerate}[(a)]
        \item ${\gamma}$ is a countable ordinal, 
        \item $\mc B=\{B_i:i<{\xi}\}$ is a base-enumeration of a
        topology ${\tau}$ on ${\gamma}$ for some  ${\xi}<{\omega}_1$,
    \item $\mc F\subs  {\gamma}\times {[\mc P({\gamma})]}^{{\omega}}$ and $|\mc F|\le {\omega}$,
    \item  If $\<{\alpha},\mc H\>\in \mc F $ 
    then every $H\in \mc H$ is ${\tau}$-dense-in-itself, and ${\alpha}$ is an accumulation point 
    of $\mc H$ in ${\tau}$.   
    \end{enumerate}
    Write $X=\<{\gamma},{\tau}\>$.
    
    If $\mbb T^{\nu}$ is a nice triple, we use the notation 
    ${\gamma}^{\nu}$, ${\tau}^{\nu}$,  $X^{\nu}$, ${\xi}^{\nu}$, $B^{\nu}_i$ for 
    $i<{\xi}^{\nu}$,  and  $\mc F^{\nu}$.

    If $\mbb T^0$ and $\mbb T^1$ are nice triples, write $\mbb T^0\preceq T^1$ iff
    \begin{enumerate}[(A)]
    \item ${\gamma}^0\le {\gamma}^1$, ${\xi}^0\le {\xi}^1$,
    \item $B^0_i=B^1_i\cap X^0$ for $i<{\xi}^0$,
    \item $B^0_i\subs B^0_k$ iff  
     $B^1_i\subs B^1_k$ for $i,k<{\xi}^0$,  
    \item $B^0_i\cap B^0_k=\empt$ iff  
     $B^1_i\cap B^1_k=\empt$ for $i,k<{\xi}^0$,    
    \item $\mc F^0\subs \mc F^1$.
    \end{enumerate}

    \begin{lemma}\label{lm:1p}
        If $\mbb T^0$ is a nice triple, $\mc H\in [\mc P(X_0)]^{\omega}$, and 
        every $H\in \mc H$ is ${\tau}_0$-dense-in-itself, then 
        there is a nice triple $\mbb T^1$ such that $\mbb T^0\preceq \mbb T^1$ and  
        $\<j,\mc H\>\in \mc F^1$ for some $j\in {\gamma}^1$.  
        \end{lemma}

 \begin{proof}
 If $\mc H$ is not locally finite in a point $x$, then 
 $\mbb T^1=\<{\gamma}^0,\mc B^0,\mc F^0\cup\{\<x,\mc H\>\}\>$ works.
    
So we can assume that $\mc H=\{H_n:n<{\omega}\}$ is locally finite. 
Hende, for each ${\gamma}^0$, we can pick a clopen  set $W_y\ni y$
such that $\{\ell: H_\ell\cap W_n\ne \empt\}$ is finite.  
Fix an enumeration $\{y_n:n<{\omega}\}$ of ${\gamma}^0$.

Thus, by induction on $n$, we can pick pairwise disjoint clopen sets 
$\{U_n:n<{\omega}\}$ and distinct  natural numbers $\ell_n$
such that 
\begin{enumerate}[(a)]
\item $\{\ell:U_n\cap H_\ell\ne \empt\}$ is finite, 
\item \label{en:t2} $U_n\cap \bigcup\{W_{y_{k}}:k\le n\}=\empt$,    
\item $U_n\cap H_{\ell_n}\ne \empt$.
\end{enumerate}
Pick $x_n\in U_n\cap H_{\ell_n}$ for $n<{\omega}$.

Let $\{B_m:m<{\omega}\}$ be an enumeration of $\mc B^0$.
By induction on $m<{\omega}$, pick 
$I_m\in {[{\omega}]}^{{\omega}}$ and $j_m<2$ such that 
\begin{enumerate}[(i)]
\item ${\omega}=I_0\supset I_1\supset\dots$,
\item $j_m=1$ iff  $\forall n\in I_{m+1} $ $x_n\in B_m$,
\item $j_m=0$ iff  $\forall n\in I_{m+1} $ $x_n\notin B_m$.
\end{enumerate}

Let $n_m=\min(I_m\setm \{n_{m'}:m'<m\})$.
Then 
\begin{displaymath}
x_{n_m}\in \bigcap\{B_{m'}:m'\le m, j_{m'}=1\}\setm \bigcup\{B_{m'}:m'\le m, j_{m'}=0\}. 
\end{displaymath}

So we can find an open $V_m$ such that 
\begin{enumerate}[(1)]
\item $x_{n_m}\in V_m\subs U_{n_m}$,
\item if $m'\le m $ and  $j_{m'}=1$, then $V_m\subs B_m$,
\item if $m'\le m $ and  $j_{m'}=0$, then $V_m\cap B_{m'}=\empt$.
\end{enumerate}

Let ${\gamma}^1={\gamma}^0+1$ and 
 ${\xi}^1={\xi}^0+{\omega}$.
For $i<{\xi}^0$, let 
\begin{displaymath}
{B^1_i}=\left\{\begin{array}{ll}
{B^0_i\cup\{{\gamma}^0\}}&\text{if $\forall^\infty m\ V_m\subs B^0_i$},\\\\
{B^0_i}&\text{if $\forall^\infty  m\  V_m\cap B^0_i=\empt$.}
\end{array}\right.
\end{displaymath} 
For ${\xi}^0\le {\eta}={\xi}^0+k<{\xi}^1$, 
let 
\begin{displaymath}
B^1_{{\xi}^0+k}=\{{\gamma}^0\}\cup\bigcup_{k<m}V_m.
\end{displaymath}
Write $\mc B^1=\{B^1_i:i<{\xi}^1\}$.
Then $\mbb T^1=\<{\gamma}^1,\mc B^1,\mc F^0\cup\{\<{\gamma}^0,\mc H\>\}\>$ meets the requirements. 

Indeed, the topology ${\tau}^1$ is $T_2$, because if $y_k\in {\gamma}^0$, 
then $W_{y_k}$ and $B^1_{{\xi}+k}$ are disjoint neighborhoods of $y_k$
and ${\gamma}^0$, respectively. 

\end{proof}       
        

    \begin{lemma}\label{lm:2p}
        If $\mbb T^0$ is a nice triple and   $A\in \dioo({\tau}^0)$, 
        $a\in {\gamma}^0\setm A$ then 
        there is a nice triple $\mbb T^1$ such that  $\mbb T^0\preceq \mbb T^1$ and 
        $a\notin \overline{A}^{{\tau}^1}$. 
        \end{lemma}
\begin{proof}
Let 
\begin{displaymath}
\mc M=\mc B^0\cup\{B^0_i\cap H: i<{\xi}^0, H\in \mc H\text{ for some } \<{\xi},\mc H\>\in \mc F, H\cap B^0_i\ne \empt\}.
\end{displaymath}
Every element of $M$ is ${\tau}^0$-dense-in-itself, and so infinite. Thus there is a  partition 
$\<C_0,C_1\>$ of ${\gamma}^0\setm (A\cup \{a\})$ such that 
$M\cap C_i\ne \empt $ for each $M\in \mc M$. 

Let ${\xi}^1={\xi}^0\dotplus {\xi}^0 \dotplus {\xi}^0$, $B^1_i=B^0_i$ for $i<{\xi}$,
$B^1_{{\xi}\dotplus i}=(B^0_i\cap C_0)\cup\{a\})$ and 
$B^1_{{\xi}\dotplus {\xi}\dotplus i}=(B^0_i\cap C_1)\cup  A$.
Then the nice  triple $\mbb T^1=\<{\gamma}^0,\mc B^1,\mc F^0\>$ satisfies the requirements.  
\end{proof}        

We will define an increasing chain $\<\mbb T^{\nu}:{\nu}\le {\omega}_1\>$ of nice triples.

Let $\{t_{\nu}:{\nu}<{\omega}_1\}$ be an  
${\omega}_1$-abundant enumeration of 
$({\omega}_1\times[{\omega}_1]^{\omega})\cup \left[{{[{\omega}_1]}^{{\omega}}}\right]$.

Let ${\gamma}^0={\omega}$, and 
let 
$\mc B^0$ be an enumeration of a  0-dimensional base of a 0-dimensional second countable 
crowded topology on ${\omega}$.
 Let $\mbb T^0=\<{\omega},\mc C^0,\empt\>$.

Assume that for some limit ordinal ${\mu}\le {\omega}_1$, $\{\mbb T^{\nu}:{\nu}<{\mu}\}$
is constructed. Let $\mc B^{\mu}=\lim_{{\nu}\to {\mu}}\mc B^{\nu}$, 
and $\mc F^{\mu}=\bigcup\{\mc F^{\nu}:{\nu}<{\mu}\}$, 
${\gamma}^{\mu}=\bigcup\{{\gamma}^{\nu}:{\nu}<{\mu}\}$.

If ${\mu}<{\omega}_1$, then the nice triple $\mbb T^{\mu}=\<X^{\mu},\mc B^{\mu},\mc F^{\mu}\>$
meets the requirements.

Assume that ${\mu}={\nu}+1$ and $\mbb T^{\nu}$ is constructed.

\smallskip
\noindent {\bf Case 1.} {$t_{\nu}=\<a,H\>\in {\omega}_1\times {[{\omega}_1]}^{{\omega}}$}.

If $a \in {\gamma}^{\nu}\setm H$,  and  $H\subs {\gamma}^{\nu}$  is ${\tau}^{\nu}$-discrete, then apply 
Lemma \ref{lm:2p} to obtain $\mbb T^{\mu}$ from $\mbb T^{\mu}$
such that 
$a\notin \overline{H}^{{\tau}_{\mu}}$.

Otherwise, let $\mbb T^{\mu}=\mbb T^{\nu}$.
\smallskip

\noindent {\bf Case 2.} $t_{\nu}\in \left[{{[{\omega}_1]}^{{\omega}}}\right]$.

If every $H\in t_{\nu}$ is a dense-in-itself  subspace of $X^{\nu}$, then 
apply 
Lemma \ref{lm:1p} to obtain $\mbb T^{\mu}$ from $\mbb T^{\mu}$
such that 
$\<a,t_{\nu}\>\in \mc F^{\mu}$ for some $a\in {\gamma}^{\mu}$.

Otherwise, let $\mbb T^{\mu}=\mbb T^{\nu}$.

\smallskip

Finally, $X^{{\omega}_1}=\<{\omega}_1,{\tau}^{{\omega}_1}\>$ satisfies the requirements of the theorem. 

Indeed, if $\mc H$ is a countable family of 
${\tau}^{{\omega}_1}$-crowded, countable subsets, then 
there is ${\nu}$ such that $t_{\nu}=\mc H$ and $\bigcup\mc H\subs {\gamma}^{\nu}$.
Thus, $\<y,\mc H\>\in \mc F^{\nu+1}$, and so $y$ is an accumulation point of 
$\mc H$ in the topology ${\tau}^{{\omega}_1}$. Hence, ${\tau}^{{\omega}_1}$
is $\nwdoo$-SP.

\end{proof}

\end{document}